\documentclass[a4paper,11pt]{article}
\usepackage{graphicx}
\usepackage{amssymb}
\usepackage{latexsym, bm}
\usepackage{multicol}
\usepackage{indentfirst}
\usepackage{amssymb,amsfonts}
\usepackage{amsmath}
\usepackage{setspace}
\newtheorem{thm}{Theorem}
\newtheorem{lem}{Lemma}

\newtheorem{fact}{Fact}

\usepackage{amssymb}
\begin{spacing}{1.01}
\begin{document}
\title{Large book--cycle Ramsey numbers}

\date{}

\author{{Qizhong Lin\footnote{Center for Discrete Mathematics, Fuzhou University,
Fuzhou, 350108 P.~R.~China. Email: {\tt linqizhong@fzu.edu.cn}. Supported in part by NSFC(No.\ 11671088) and ``New century excellent talents support plan for institutions of higher learning in Fujian province"(SJ2017-29).}
~and
 {Xing Peng\footnote{School of Mathematical Sciences, Anhui University, Hefei 230601, P.~R.~China. E-mail: {\tt x2peng@ahu.edu.cn}. Supported by the NSFC grants
(No.\ 11601380, 12071002) and a Start--up Fund from Anhui University.}}}}
\maketitle

\maketitle
\begin{abstract}
Let $B_n^{(k)}$ be the book graph which consists of $n$ copies of $K_{k+1}$  all  sharing a common $K_k$, and let $C_m$ be a cycle of length $m$. In this paper, we first determine the exact value of $r(B_n^{(2)}, C_m)$ for $\frac{8}{9}n+112\le m\le \lceil\frac{3n}{2}\rceil+1$ and $n \geq 1000$.
This answers a question of Faudree, Rousseau and Sheehan  (Cycle--book Ramsey numbers, {\it Ars Combin.,} {\bf 31} (1991), 239--248) in a stronger form when $m$ and $n$ are large.
  Building upon this exact result, we are able to determine the asymptotic value of $r(B_n^{(k)}, C_n)$ for each  $k \geq 3$. Namely, we prove that for each $k \geq 3$, $r(B_n^{(k)}, C_n)= (k+1+o_k(1))n.$
This extends a result due to Rousseau and Sheehan (A class of Ramsey problems involving trees, {\it J.~London Math.~Soc.,} {\bf 18}  (1978), 392--396).

\medskip




\end{abstract}

\section{Introduction}\label{sec-1}
For graphs $H_1$ and $H_2$, the Ramsey number $r(H_1, H_2)$ is the minimum integer $N$ such that every red-blue edge coloring of the complete graph $K_N$ contains either a red $H_1$ or a blue $H_2$.
Let $B_n^{(k)}$ be the book graph which consists of $n$ copies of $K_{k+1}$  all  sharing a common $K_k$. When $k=2$, we write $B_n$ instead of $B_n^{(2)}$ for convenience.
Book Ramsey numbers have attracted a lot of  attention, see \cite{efrs1,rs,tho,nr1,nr2,nrs} and other related references.
In particular, answering a  question of Erd\H{o}s et al. \cite{efrs1}, Conlon \cite{con} established an asymptotic version of Thomason's conjecture \cite{tho} by showing
\[
r(B_n^{(k)},B_n^{(k)})=(2^k+o_k(1))n.
\]
The upper bound was improved to $2^{k} n+O_k\big(\frac{n}{(\log\log\log n)^{1/25}}\big)$ by using a different method, see Conlon, Fox and Wigderson \cite{cfw}.

Let $C_m$ and $T_m$ be a cycle and a tree of order $m$, respectively.
The Ramsey numbers of book versus tree and book versus cycle also received a great deal of attention.
Strengthening a classical result due to Chv\'{a}tal \cite{chv}, Rousseau and Sheehan \cite{rs2}  established that
\begin{align}\label{r-s}
r(B_n^{(k)}, T_n)= (k+1)(n-1)+1.
\end{align}
For more book-tree Ramsey numbers, see e.g. \cite{rs2,efrs,bef,frs1}.
A natural question is  whether we can prove a similar result  for book-cycle Ramsey number when they have nearly equal order.

The study of book-cycle Ramsey numbers goes back to \cite{rs} by Rousseau and Sheehan. In particular, they proved $r(B_n,C_3)=2n+3$ for $n>1$.
In \cite{frs0,frs}, Faudree, Rousseau and Sheehan  proved some results for $r(B_n,C_4)$, and generally,
 \[
r(B_n,C_m) = \left\{ \begin{array}{ll}
2n+3 & \textrm{if $m\ge5$ is  odd and $m\le \frac{n}{4}+\frac{13}{4}$,}\\
2m-1 & \textrm{if $m\ge2n+2$.}
\end{array} \right.
\]
 Improving upon the result in \cite{frs}, Shi \cite{shi} obtained that $r(B_n,C_m)=2m-1$ for $m>\frac{3n}{2}+\frac{7}{4}$.
 In the same paper,  the author also obtained  $r(B^{(3)}_n,C_m)=3m-2$ for $m>\max\{\frac{6n+7}{4},70\}$.
For fixed $k\ge1$ and odd $m\ge 3$, Liu and Li \cite{ll} proved that $r(B^{(k)}_n,C_m)=2(n+k-1)+1$ when $n$ is large.
   One can easily see that $r(B_n, C_m)>3(n-1)\ge \max\{2m-1,2n+3\}$ for $6\le n\le m\le\frac{3n}{2}-1$,
  and $r(B_n, C_m)>3(m-1)\ge \max\{2m-1,2n+3\}$ for $\frac{2n}{3}+2\le m\le n$.
This suggests that the formula for $r(B_n,C_m)$ varies when $m$ and $n$ change, especially when $m$ and $n$ are nearly equal.
 As mentioned in \cite{frs},
 ``the problem of computing $r(B_n,C_m)$ when $m$ is odd and $m$ and $n$ are nearly equal provides an unanswered test of strength''.

 The goal of this paper is to study the Ramsey number $r(B_n^{(k)},C_m)$ when  $n$ and $m$ are nearly equal.
First, we determine the exact value of $r(B_n,C_m)$ for $\frac{8}{9}n+112 \le m\le \lceil\frac{3n}{2}\rceil+1$ and $n \geq 1000$,
which provides an answer to the question by Faudree, Rousseau and Sheehan \cite{frs} in a stronger form when $m$ and $n$ are large.

\begin{thm}\label{Bn-Cm}
For  $n \geq 1000$,
\[
r(B_n,C_m) = \left\{ \begin{array}{ll}
3m-2 & \textrm{if $\frac{8n}{9}+112\le m\le n$,}\\
3n-1 & \textrm{if $m=n+1$,}\\
3n & \textrm{if $n+2\le m\le \frac{3n+1}{2}$,}
\\
2m-1 & \textrm{if $m= \lceil\frac{3n}{2}\rceil+1$.}
\end{array} \right.
\]
\end{thm}

\noindent
{\bf Remark 1.} We observe that the formula for $r(B_n,C_m)$ undergos phase transitions when $m \in \{n,n+1,n+2\}$.

So far, the value of $r(B_n,C_m)$ is known for $m \geq \frac{8n}{9}n+112$
as well as odd $m$ with $m\le \frac{n}{4}+\frac{13}{4}$.
It requires new ideas to determine the value of $r(B_n,C_m)$ when $m$ and $n$ are in other ranges.

Based on Theorem \ref{Bn-Cm}, we extend \eqref{r-s} by showing the asymptotic value of $r(B_n^{(k)}, C_n)$ for each fixed integer $k \geq 3$ as follows.

\begin{thm}\label{main}
Let $k\ge 3$ be a fixed integer. We have $$r(B_n^{(k)}, C_n)= (k+1+o_k(1))n.$$
\end{thm}
It is a challenge to determine the exact value of $r(B_n^{(k)},C_m)$ when $m$ and $n$ are nearly equal for each $k \geq 3$.

Throughout this paper, we will use the following notation.
Let  $G$ be a graph with vertex set $V$.  For each vertex $v\in V$, we use $N_G(v)$ and $d_G(v)$ to denote the {\it neighborhood} and the {\it degree}, respectively.
If $U\subseteq V$, then  $G[U]$  denotes the subgraph induced by $U$. Moreover,
let  $N_G(v,U)=N_G(v)\cap U$ and $d_G(v,U)=|N_G(v,U)|$.
 The graph $G-v$ is the one obtained from $G$ by deleting the vertex $v$ and all edges incident to $v$.

The rest of the paper is organized as follows. In Section 2, we will collect several results  which will be used to prove our results.  In Section 3, we will present proofs of Theorem \ref{Bn-Cm} and Theorem  \ref{main}.

\section{Preliminaries}
In this section, we collect a number of previous results which are needed for our proofs. Crucial tools include a refined version of Szemer\'edi's regularity lemma and the (weakly) pancyclic properties of graphs.

\subsection{Regularity method}\label{reg-lem}
Szemer\'edi's regularity lemma \cite{sze78} is a powerful tool in extremal graph theory.  The regularity lemma is also called the uniformity lemma, see e.g., Bollob\'as \cite{bol98} and Gowers \cite{gow}.
 Applications of the regularity method are fruitful.
 We refer the reader to the survey of Koml\'os and Simonovits \cite{kom} and other related references.

A key tool in the proof of Theorem \ref{main} is a refined version of the Szemer\'edi's regularity lemma  by Conlon \cite{con}. We state it and its related results as following.

Let $G$ be a graph defined on vertex set $V=V(G)$.
 For $X, Y\subseteq V$, denote $e_G(X,Y)$ by the number of pairs in $X \times Y$ that are edges of $G$. The ratio \label{'edge-b-d'}
$$
d_G(X,Y) = \frac{e_G(X,Y)}{|X||Y|}
$$
is called the {\it edge density} of $(X,Y)$ in $G$, which can be understood as the probability that a random pair $(x,y)$  from $X\times Y$ is an edge. If $X\cap Y\neq \emptyset$, then  edges in $X\cap Y$ are counted twice.

For $\epsilon >0$, a pair $(U,W)$ of nonempty  sets
$U,W\subseteq V$ is called {\it $\epsilon$-regular} if
\[
|d_G(X,Y)-d_G(U,W)|\le \epsilon
\]
for every $X\subseteq U, Y\subseteq W$ such that $|X|\ge \epsilon |U|$ and $|Y|\ge \epsilon |W|$. We say a subset $U$ is $\epsilon$-regular if the pair $(U,U)$ is $\epsilon$-regular.


An  {\em equitable partition} of a graph $G$ is a partition $V (G) = \sqcup^m_{i=1}V_i$ of the vertex set of $G$ such that $\big||V_i|-|V_j|\big|\le 1$ for all $i$ and $j$.

\medskip
We mention two properties for $\epsilon$-regular pairs, see e.g. \cite{kom}.
\begin{fact}\label{lar-dg}
 Let $(U,W)$ be an $\epsilon$-regular pair with edge density $d$. If  $Y\subseteq W$ with $|Y|\ge \epsilon |W|$, then there
exists a subset $U'\subseteq U$ with $|U'|\ge (1-\epsilon)|U|$ such that each vertex in $U'$ is adjacent to at least $(d-\epsilon)|Y|$
vertices in $Y$.
\end{fact}

 \begin{fact}\label{hdt}
 Let $(U,W)$ be an $\epsilon$-regular pair in graph $G$. If $X\subseteq U$, $Y\subseteq W$ with $|X|\geq \gamma|U|$ and $|Y|\geq \gamma|W|$ for some $\gamma>\epsilon$, then $(X,Y)$ is $\epsilon'$-regular such that $|d_G(U,W)-d_G(X,Y)|\le\epsilon$, where $\epsilon'=\max\{\epsilon/\gamma, 2\epsilon\}$.
\end{fact}

We need the following refined version of the regularity lemma by Conlon \cite[Lemma 3]{con}.
In the same spirit, to prove an induced removal lemma,  Alon et.~al \cite{afk} obtained a result in which all pairs $(W_i,W_j)$ are $\epsilon$-regular.
 \begin{lem} \label{con0}
For every $0<\epsilon<1$ and natural number $m_0$, there exists a natural number $M$ such that every graph $G$ with at least $m_0$ vertices has an equitable partition $V(G)=\sqcup_{i=1}^m V_i$ with $m_0 \leq m \leq M$ parts and subsets $W_i \subset V_i$ such that $W_i$ is $\epsilon$-regular for all $i$ and, for all but $\epsilon m^2$ pairs $(i,j)$ with $1\leq i \neq j \leq m, (V_i,V_j)$, $(W_i, V_j)$ and $(W_i,W_j)$ are $\epsilon$-regular with $|d_G(W_i,V_j)-d_G(V_i,V_j)| \leq \epsilon$ and $|d_G(W_i,W_j)-d_G(V_i,V_j)| \leq \epsilon$.
\end{lem}

We will also use the following counting lemma from  \cite[Lemma 5]{con}.
\begin{lem}\label{con}
For any $\delta>0$ and any natural number $k$, there is $\eta>0$ such that if $U_1,\dots,U_k$, $U_{k+1},\dots,U_{k+\ell}$ are (not necessarily distinct) vertex sets with $(U_i,U_{i'})$ $\eta$-regular of density $d_{i,i'}$ for all $1\le i< i'\le k$ and $1\le i\le k< i'\le k+\ell$ and $d_{i,i'}\ge\delta$ for all $1\le i< i'\le k$, then there is a copy of $K_k$ with vertex $u_i\in U_i$ for each $1\le i\le k$ which is contained in at least
$$
\sum_{j=1}^\ell\left(\prod_{i=1}^kd_{i,k+j}-\delta\right) |U_{k+j}|
$$
labeled copies of $K_{k+1}$ with vertex $u_{k+1}$ in $\cup_{j=1}^\ell U_{k+j}$.
\end{lem}

The next lemma due to Benevides and Skokan \cite{Benevides-2009} is a stronger version of the original one by {\L}uczak \cite[Claim 3]{Luczak-1999}. Both have similar proofs by using Fact \ref{lar-dg}.

\begin{lem}\label{ben-s}
For every $0<\beta< 1$, there exists an $n_0$ such that for every $n > n_0$ the following holds: Let $G$ be a bipartite graph with bipartition $V(G)=U \cup W$ such that $|U|=|W|=n$. Furthermore, let the pair $(U,W)$ be $\epsilon$-regular with density at least $\beta$ for some $\epsilon$ satisfying $0<\epsilon<\beta/100$. Then for each $\ell$, $1\le \ell \le n-5\epsilon n/\beta$, and for each pair of vertices $u \in U$, $w \in W$ with $d_G(u)\ge 4\beta n/5$ and $d_G(w)\ge 4\beta n/5$, $G$ contains a path of length $2\ell+1$ connecting $u$ and $w$.
\end{lem}

 The following lemma will be used to find long odd cycles in graphs.
\begin{lem}\label{edge}
Suppose that $(A,B)$ and  $(B,C)$  are $\epsilon$-regular pairs with density at least $\beta$ in a graph $G$,
here $\epsilon<\beta/50$  and $A,B$ and $C$ are pairwise disjoint.
 If $|A|=|C|=n$ and $|E(B)| \geq |B|^2/5$,
 then there is an edge $xy \in E(B)$ such that $d(x,A) \geq (\beta-\epsilon)n$ and $d(y,C) \geq (\beta-\epsilon)n$.
\end{lem}
\noindent
{\bf Proof:} Let $H$ be the $2\epsilon |B|$-core of  the subgraph induced by $B$, i.e., $H$ is the maximum induced subgraph of $B$ with minimum degree at least $2\epsilon |B|$. As $|E(B)| \geq |B|^2/5$ and $\epsilon$ is small enough, we can see $H$ is not empty and $|V(H)| \geq 2\epsilon |B|$. Since $(A,B)$ is an $\epsilon$-regular pair, all but at most $\epsilon |B|$ vertices in $H$ have at least $(\beta-\epsilon) n$ neighbors in $A$ from Fact \ref{lar-dg}. We assume $x$ is such a vertex. The definition of $H$ gives $N_H(x) \geq 2\epsilon |B|$. Now, $(B,C)$ being an $\epsilon$-regular pair yields that there is a vertex $y \in N_H(x)$ having at least $(\beta-\epsilon) n$ neighbors in $C$. The edge $xy$ is a desired one and the proof is complete. \hfill$\Box$

\medskip
The following lemma by \L ucazk \cite[Claim 7]{Luczak-1999} is a key ingredient in the proof of Theorem \ref{main}.
\begin{lem} \label{tluczak}
For every $0<\delta<10^{-15}, \alpha > 2\delta$ and $n \geq \rm{exp}(\delta^{-16}/\alpha)$ the following holds. Each graph $G$ on $n$ vertices which contains no odd cycles longer than $\alpha n$ contains subgraphs $G'$ and $G''$ such that:

\smallskip

$(1)$ $V(G') \cup V(G'')=V(G)$, $V(G') \cap V(G'')=\emptyset$ and each of the sets $V(G')$ and $V(G'')$ is either empty or contains at least $\alpha\delta n/2$ vertices;

\smallskip

$(2)$ $G'$ is bipartite;

\smallskip

$(3)$ $G''$ contains not more than $\alpha n |V(G'')|/2$ edges;

\smallskip

$(4)$ all but at most  $\delta n^2$ edges of $G$ belong to either $G'$ or $G''$.

\end{lem}

\subsection{Pancyclic properties of graphs}\label{pan-cyc}

For a graph $G$, we use $g(G)$ and $c(G)$ to denote its \emph{girth} and \emph{circumference}, i.e., the length of a shortest cycle and a longest cycle of $G$. Similarly, the {\it odd girth} of $G$ is the length of a shortest odd cycle in $G$.
 A graph is called \emph{weakly pancyclic} if it contains cycles of every length between its  girth and its circumference. A graph is \emph{pancyclic} if it is weakly pancyclic with girth $3$ and circumference $n=|V(G)|$. We say a graph is {\it 2-connected} if it remains connected after the deletion of any  vertex.

For a graph $G$, let $\delta(G)$ denote the minimum degree of $G$. The following classical result is due to Dirac \cite{di}.
\begin{lem}\label{1-thm-3}
Let $G$ be a 2-connected graph of order $n$ with minimum degree $\delta=\delta(G)$. Then $c(G)\geq \min\{2\delta,n\}$.
\end{lem}

Dirac's result tells us that the circumference of a 2-connected graph cannot be too small.
In particular, if $\delta=\delta(G)\ge n/2$, then $c(G)=n$. This is a well-known result for a graph being hamiltonian.
For the special case of $\delta \geq n/2$,
the following result due to Bondy \cite{Bondy} tells us more about the structure of the graph.

\begin{lem}\label{1-thm-1}
If a graph $G$ with $n$ vertices satisfies  $\delta(G)\geq n/2$, then $G$ is pancyclic unless $n=2r$ and $G=K_{r,r}$.
\end{lem}

The following is an elegant extension on graphs being weakly pancyclic by Brandt, Faudree and Goddard \cite{Brandt}, which is a key ingredient in the proofs of Theorem \ref{Bn-Cm} and Theorem \ref{main} for  $k=3$.

\begin{lem}\label{1-thm-2}
Let $G$ be a 2-connected nonbipartite graph of order $n$ with $\delta(G)\geq n/4+250$. Then $G$ is weakly pancyclic unless $G$ has odd girth 7, in which case it has every cycle from 4 up to its circumference except the 5-cycle.
\end{lem}
We will also need the following simple fact which can be seen using the Breadth-First-Search.
\begin{fact}\label{f1}
If a graph $G$ with $n$ vertices satisfies $\delta(G)\geq cn$ for some constant $c>0$, then $g(G) \leq 4$ provided $n>c^{-2}$.
\end{fact}

\section{Proofs of Theorems \ref{Bn-Cm} and  \ref{main}}

In this section, we will give  proofs for our main results.
Throughout the proof, when considering a red-blue edge coloring of $K_N$, we always use $R$ and $B$ to denote  subgraphs formed by red and blue edges, respectively. We also suppose that $n \geq 1000$ for Theorem \ref{Bn-Cm} and $n$ is sufficiently large for Theorem \ref{main}.

\subsection{Proof of Theorem \ref{Bn-Cm}}

We first give the following simple fact.
\begin{fact}\label{bs}
Let $G$ be a graph which consists of three connected components $V_1,V_2$ and $V_3$.
\smallskip

(i) If the largest connected component has at least $n$ vertices, then the complement $\overline{G}$ contains a $B_n$.

(ii) If one of sets  $V_1,V_2$ and $V_3$ contains a non-edge while the other two sets have at least $n$ vertices in total, then the complement $\overline{G}$ contains a $B_n$.
\end{fact}

The proof of Theorem \ref{Bn-Cm} contains three parts.

\bigskip
\noindent
{\bf (I) $\frac{8}{9}n+112 \le m\le n$}
\medskip

The lower bound  $r(B_n, C_m)>3m-3$ holds since the graph with three disjoint copies of $K_{m-1}$ contains no $C_m$ and its complement contains no $B_n$. Thus it suffices to prove the upper bound.
Let $N=3m-2$, and consider a red-blue edge coloring of $K_N$  on vertex set $V$.

\medskip
\noindent
{\bf Case 1.} There is a vertex $v \in V$ with $d_R(v)\geq 2n$.
 \medskip

 We choose a subset $X\subseteq N_R(v)$ with $|X|=2n$.
If  there is a vertex $x\in X$ with $d_R(x,X)\geq n$, then the red subgraph induced by $x$ and $v$ together with their common neighbors in $X$ contains a  $B_n$. Thus we assume  $d_B(x,X)\geq n$ for each $x \in X$, i.e. $\delta(B[X])\ge n$.
By Lemma \ref{1-thm-1}, $B[X]$ is pancyclic or $B[X]=K_{n,n}$.
There will be a blue $C_m$ if $B[X]$ is pancyclic, so we assume $B[X]=K_{n,n}$ with color classes $X_1$ and $X_2$.
If there exists a vertex $y\in V\setminus (X\cup \{v\})$ such that $d_B(y,X_i)\ge1$ for each $i \in \{1,2\}$, then $B[X\cup \{y\}]$ contains blue cycles of length between 3 and $2n+1$ and it definitely contains a blue $C_m$. Thus  each vertex of $V\setminus (X\cup \{v\})$ is completely red-adjacent to $X_1$ or $X_2$. Suppose that $z \in V\setminus (X\cup \{v\})$ is red-adjacent to $X_1$. It follows that $R[X_1\cup \{z,v\}]$ contains a red $K_{n+2}-e$ and definitely a red $B_n$. We are through in this case.

\medskip
\noindent
{\bf Case 2.} For each vertex $v \in V$, $d_R(v) \leq  2n-1$. i.e., $\delta(B)\geq 3m-2n-2$.
\medskip

 If $B$ is bipartite, then the larger color class of $B$ contains at least $\frac{3m-2}{2}\geq n+2$ vertices and induces a red clique of size at least $n+2$. Therefore, there is a red $B_n$.

If $B$ is 2-connected, then we can find a blue $C_m$ as follows. Since
$$\delta(B)\geq 3m-2n-2\geq \frac{3m-2}{4}+250$$
provided $m\ge(8n+1006)/9$ which is guaranteed by the assumption that $m\ge \frac{8n}{9}+112$,
it follows from Lemma \ref{1-thm-2} that $B$ is weakly pancyclic unless $B$ has odd girth 7, in which case it contains every cycle of length from 4 up to its circumference $c(B)$ except the 5-cycle.
 Moreover, by Lemma \ref{1-thm-3}, $c(B)\geq 2(3m-2n-2)\ge m$.
Note that Fact \ref{f1} implies $g(B) \leq 4$ since $\delta(B)\geq 3m-2n-2>N/4$.
 Thus  there is  a blue $C_m$.

 In the following, we assume that $B$ is nonbipartite and not 2-connected.  Suppose that $B-u$ is disconnected for some vertex $u \in V$, here it includes the case where $B$ is disconnected. Since $\delta(B-u)\geq 3m-2n-3$, we have that each connected component has at least $3m-2n-2$ vertices. So there are at most three connected components in $B-u$. Otherwise, $4(3m-2n-2)>3m-3=(N-1)$, which is a contradiction.

\medskip
 {\textbf{Subcase 2.1.}}  $B-u$ contains three connected components.
 \medskip

Let $V_1,V_2$ and $V_3$ be the vertex sets of these three connected components of $B-u$. We assume that $V_3$ is the largest one. If $|V_3| \geq m$, then we can find a blue $C_m$ as follows. Since each connected component has size at least $3m-2n-2$,
 we have
\[
|V_3|\le (3m-3)-2(3m-2n-2)=4n-3m+1.
\]
Since $\delta(B[V_3])\ge \delta(B)-1 \geq 3m-2n-3>|V_3|/2$, Lemma \ref{1-thm-1} implies that $B[V_3]$ is pancyclic and contains a blue $C_m$. Thus we assume $|V_3| \leq m-1$. As  $V_3$ is the largest connected component and $|V\setminus \{u\}|=3m-3$, we get  $|V_i|= m-1$ for each $1\le i \le3$.

We claim that each $V_i$ induces a blue clique $K_{m-1}$.
Otherwise, $R[V_1\cup V_2\cup V_3]$ contains a $B_n$ by Fact \ref{bs}(ii) since $2m-2\ge n$.
Because $d_B(u)\ge 3m-2n-2\ge4$, we get that $u$ has at least two blue neighbors in $V_i$ for some $1 \leq i \leq 3$, say $V_1$.
 Therefore, $B[V_1 \cup\{u\}]$ contains a $C_m$.

\medskip
{\textbf{Subcase 2.2.}}  $B-u$ contains exactly two connected components.
\medskip

Let $V_1$ and $V_2$ be the vertex sets of these two connected components with $|V_1|\le |V_2|$.
Clearly, $|V_1|\le \frac{3m-3}{2}$.
If $|V_1| \geq m$, then $B[V_1]$ contains a  $C_m$.
Indeed, Lemma \ref{1-thm-1} implies that $B[V_1]$ is pancyclic as $\delta(B[V_1])\ge 3m-2n-3>|V_1|/2$ .
Hence, we assume $3m-2n-2\le|V_1|\le m-1$ as each connected component has at least $3m-2n-2$ vertices.
Clearly, $2m-2\le |V_2|\le 2n-1$.

If $B[V_2]$ is bipartite with color classes $X$ and $Y$ satisfying $|X|\ge|Y|$, then we can find a red $B_n$  as following.  We notice
$$|V_1|+|X|\ge |V_1|+\frac{|V_2|}{2}>\frac{|V_1|+|V_2|}{2}\ge \frac{3m-3}{2}\ge n+2.$$
Since $X$ induces a red clique with $|X|\ge2$ and all edges between $V_1$ and $X$ are red,
$R[V_1\cup X]$ contains a  $B_n$.

If $B[V_2]$ is 2-connected, then we have
$\delta(B[V_2])\ge 3m-2n-3\ge\frac{2n-1}{4}+250\ge\frac{|V_2|}{4}+250$
by noting that $m\ge \frac{8n}{9}+112$ and $n \geq 1000$.
Lemma \ref{1-thm-3}, Lemma \ref{1-thm-2} and Fact \ref{f1} imply that there is  a blue $C_m$ in $B[V_2]$.

Therefore, we are left to consider the case where $B[V_2]$ is nonbipartite and contains a cut vertex. Suppose that $B[V_2]-w$ is disconnected for some vertex $w \in V_2$,  here it includes the case where $B[V_2]$ is not connected. As $\delta(B[V_2])\ge 3m-2n-3$, it follows that $B[V_2]-w$ contains exactly two connected components, denoted by $V_2'$ and $V_2''$.

Note that $|V_2'|,|V_2''|\le |V_2|-(3m-2n-3)\le 4n-3m+2$.
If either $|V_2'| \geq m$ or $|V_2''| \geq m$, say $|V_2'| \geq m$,
then Lemma \ref{1-thm-1} implies that $B[V_2']$ is pancyclic and contains a blue $C_m$ since $\delta(B[V_2'])\ge 3m-2n-4>|V_2'|/2$.
Thus we assume  $|V_2'| \leq m-1$ and $|V_2''|\le m-1$.
Note that $|V_1\cup V_2'\cup V_2''|=3m-4$ and $\max\{|V_1|,|V_2'|,|V_2''|\}\le m-1$. We get
$$m-2\le|V_1|,|V_2'|,|V_2''|\le m-1.$$

We claim that each of $V_1$, $V_2'$ and $V_2''$ induces a blue clique. Otherwise,
Fact \ref{bs}(ii) implies that  $R[V_1\cup V_2'\cup V_2'']$ contains a $B_n$
 by noting $2m-3\ge n$.

If $|V_1|=m-2$, then $|V_2'|=|V_2''|=m-1$. Since $\delta(B[V_2])\ge 3m-2n-3\ge 4$,
we have either $d_B(w,V_2')\ge2$ or $d_B(w,V_2'')\ge2$ and so either $B[V_2'\cup \{w\}]$ or $B[V_2''\cup \{w\}]$ contains a  $C_m$.

If $|V_1|=m-1$, then we can assume $|V_2'|=m-1$ and $|V_2''|=m-2$ without loss of generality.
 If $d_B(u,V_1)\ge2$, then $B[V_1\cup \{u\}]$ contains a $C_m$. Hence $d_B(u,V_1)\le1$.
 Similarly, we have $d_B(u,V_2')\le1$. Thus $d_B(u,V_2'')\ge3$ by noting that $\delta(B)\ge3m-2n-2\ge6$.
Repeating the argument above,  we can show that $d_B(w,V_2'')\ge3$.
Now, $B[V_2''\cup\{u,w\}]$ contains a $C_m$ as desired.

We proved part I of Theorem \ref{Bn-Cm}.
\hfill$\Box$

\bigskip\noindent
{\bf (II) $m=n+1$}
\medskip

Let $G$ be the graph which consists of three $K_n$ sharing a common vertex.
The lower bound $r(B_n,C_{n+1})>3n-2$ follows from the fact that $G$ contains no $C_{n+1}$ and its complement contains no $B_n$.
  To show the  upper bound $r(B_n,C_{n+1})\le3n-1$, we consider a red-blue edge coloring of $K_{3n-1}$ on vertex set $V$.

  We follow the proof for Part I step by step. We will end up with the case which corresponds to Subcase 2.1 of the proof for Part I. In the following, we suppose $\delta(B)\ge n-1$. If there is a vertex $u\in V$ such that $B-u$ has three connected components, then we can easily find a red $B_n$  by Fact \ref{bs}(i) since the largest connected component must have order at least $\lceil(3n-2)/3\rceil=n$.

  Therefore, we suppose that there is a vertex $u$ such that $B-u$ has exactly two connected components $V_1$ and $V_2$ with $|V_1| \leq |V_2|$. Furthermore, similar to Subcase 2.2, we can assume that $B[V_2]$ is nonbipartite and there is some vertex $w \in V_2$ such that $B[V_2]-w$ is disconnected. The assumption $\delta(B-u) \geq n-2$ yields that $|V_1|,|V_2| \geq n-1$. Similarly,  $B[V_2]-w$ has two connected components, say $V_2'$ and $V_2''$, which satisfy $|V_2'|  \geq |V_2''| \geq n-2$.

If $|V_1|\ge n$, then we can definitely find a red $B_n$  by Fact \ref{bs}(i).
  Thus we assume $|V_1|=n-1$. As $|V_2'\cup V_2''|=2n-2$ and
  $|V_2'|  \geq |V_2''| \geq n-2$, we get either $|V_2'|=n$ and $|V_2''|=n-2$, or $|V_2'|=|V_2''|=n-1$. In the former case, a red $B_n$ is ensured again  by Fact \ref{bs}(i).
  In the latter case, we have that each of $V_1$, $V_2'$ and $V_2''$ induces a blue clique $K_{n-1}$. Otherwise, Fact \ref{bs}(ii) gives a red $B_n$
  by noting $2n-2\ge n$.


Since $u$ is a cut-vertex and $w \in V_2$, we get $w$ is completely red-adjacent to $V_1$. Let $x\in V_1$ be a red neighbor of $w$. If $w$ has another red neighbor $y$ in either $V_2'$ or $V_2''$, say $V_2'$, then $R[V_2''\cup\{x,y,w\}]$ contains a  $B_n$, where $xy$ is the edge shared by $n$ triangles. Hence $w$ is completely blue-adjacent to $V_2'\cup V_2''$.
 If $d_B(u,V_2')\ge2$ or $d_B(u,V_2'')\ge2$, say $V_2'$,  then $B[V_2'\cup\{u,w\}]$ contains a $C_{n+1}$. Otherwise, we take two red neighbors $a$ and $b$ of $u$, where $a\in V_2'$ and $b\in V_2''$. Clearly, $R[V_1\cup\{a,b,u\}]$ contains a $B_n$, where $ab$ is the edge shared by $n$ triangles.

 The proof of Part II is complete.

\bigskip\noindent
{\bf (III) $n+2\le m\le \frac{3n+1}{2}$}
\medskip

The lower bound $r(B_n,C_{m})>3n-1$ can be seen as follows.
Let $K_{n-1,n-1,n-1}$ be the complete tripartite graph with color classes $U_1,U_2$ and $U_3$. Let $s$ and $t$ be two new vertices.
  If $G$ is a graph obtained from $K_{n-1,n-1,n-1}$ by adding all edges between $s$ and $U_1$, and all edges between $t$ and $U_3$, then
 $G$ contains no $B_n$ and its complement contains no $C_m$. The lower bound follows.
 To show the upper bound $r(B_n,C_{m})\le3n$, we consider a red-blue edge coloring of $K_{3n}$ on vertex set $V$.

We assume $\delta(B)\geq n$, since the proof is similar to Case 1 of Part I if there is a vertex $u$ with $d_R(u)\geq 2n$.

 If the blue graph $B$ is bipartite, then the larger color class of $B$ induces a red clique of size at least $n+2$ and hence there is a red $B_n$.
 So we assume $B$ is nonbipartite.
If further the blue graph $B$ is 2-connected, then the existence of a blue $C_m$ follows from Lemma \ref{1-thm-3}, Lemma \ref{1-thm-2} and Fact \ref{f1} since $\delta(B)\geq n\geq \frac{3n}{4}+250$ for $n\ge1000$.

Therefore, we need only to consider the case where $B$ is nonbipartite and not 2-connected. Let  $u \in V$ be a vertex such that $B-u$ is disconnected. Since $\delta(B-u)\geq n-1$, each connected component of $B-u$ has at least $n$ vertices and hence $B-u$ has exactly two connected components.

Let $V_1$ and $V_2$ be the vertex sets of these two connected components with $|V_1|\le |V_2|$. If $|V_1| \geq m$, then we can find a blue $C_m$ as following.
Note that $|V_1|\le \frac{3n-1}{2}$.  We have $\delta(B[V_1])\ge n-1>|V_1|/2$ for $n\ge6$. Lemma \ref{1-thm-1} again implies that $B[V_1]$ is pancyclic, and so there is a blue $C_m$.
Thus we assume $n\le|V_1| \leq m-1$ and $3n-m\le |V_2|\le 2n-1$.

If $B[V_2]$ is bipartite with bipartition $(X,Y)$, where $|X|\ge|Y|$, then
we claim  $R[V_1 \cup X]$  contains a  $B_n$. To see this,
we notice $$|V_1|+|X|\ge |V_1|+\frac{|V_2|}{2}\ge \frac{|V_1|+|V|-1}{2} \geq \frac{4n-1}{2} \geq  n+2.$$
Since $X$ induces a red clique and all edges between $V_1$ and $X$ are red and $|X|\ge2$, there is a red $B_n$ in $R[V_1 \cup X]$ as claimed.

Moreover, if $B[V_2]$ is 2-connected, then the existence of  a blue $C_m$ in $B[V_2]$ again follows from Lemma \ref{1-thm-3}, Lemma \ref{1-thm-2} and Fact \ref{f1} since
$\delta(B[V_2])\ge n-1 \geq \frac{2n-1}{4}+250 \ge\frac{|V_2|}{4}+250$ for $n\ge502$.

 Therefore, we assume $B[V_2]$ is nonbiparite and  not 2-connected.
  Let $w\in V_2$ be a vertex such that $B[V_2]-w$ is disconnected.
   This means that the blue subgraph induced by $V\setminus\{u,w\}$ contains three connected components exactly
   since $\delta(B[V\setminus\{u,w\}])\ge n-2 > {|V\setminus\{u,w\}|}/{4}$.
Now, as the largest connected component has at least $\lceil(3n-2)/3\rceil\ge n$ vertices,
we can find a desired red $B_n$ in $V\setminus\{u,w\}$ by Fact \ref{bs}(i).

This completes the proof for part III.
\hfill$\Box$

\bigskip\noindent
{\bf (IV) $m= \lceil\frac{3n}{2}\rceil+1$}
\medskip

The lower bound  $r(B_n, C_m)>2m-2$ is clear since the graph with two disjoint $K_{m-1}$ contains no $C_m$ and its complement contains no $B_n$.
For the upper bounds $r(B_n,C_{m})\le2m-1=3n+1$ if $n$ is even and $r(B_n,C_{m})\le2m-1=3n+2$ if $n$ is odd,
one can easily follow the proof of Part III step by step apart from a few modifications.

This completes the proof for part IV and hence the proof of Theorem \ref{Bn-Cm}.
\hfill$\Box$

\subsection{Proof of Theorem \ref{main}}

We note for $k\ge3$, the graph with $k+1$ disjoint copies of $K_{n-1}$ contains no $C_n$ and its complement contains no $B_n^{(k)}$, so we have
$
r(B_n^{(k)}, C_n)\ge(k+1)(n-1)+1.
$
Therefore, it suffices to establish the upper bound.
The proof is by induction on $k\ge3$.  The base case where $k=3$ is built upon an induction idea and the case of $k=2$.

\bigskip\noindent
{\bf Step 1: $B_n^{(3)}$ versus $C_n$}

\medskip
Let $0<\xi<1/10$ and  $N=\lfloor(4+\xi)n\rfloor$. Consider a red-blue edge coloring of $K_N$ on vertex set $V$.
If there exists a vertex $v\in V$ with $d_R(v)\geq 3n-2$, then by Theorem \ref{Bn-Cm}, we can find either a blue $C_n$ or  a red $B_n$ in the red neighborhood of $v$.
We are done if there is a blue $C_n$, so we assume that there is a red $B_n$ in the red neighborhood of $v$.
Now, this red $B_n$ together with $v$ form a red $B^{(3)}_n$. 
   Thus we are left to consider the case where  $\delta(B)> (1+\xi)n.$ Fact \ref{f1} implies $g(B) \leq 4$.

We claim $B$ is nonbipartite. Otherwise, one of its color classes induces a red clique of size at least $N/2\ge n+3$, which will give us a red $B_n^{(3)}$.

Moreover, we can assume that $B$ is 2-connected. Otherwise, suppose that
there exists a vertex $u$ such that $B-u$ is disconnected. Then $B-u$ has two or three connected components as $\delta(B)> (1+\xi)n.$
 Let $V_1\subseteq V$ be the smallest connected component of $B-u$. Recall the assumption $\delta(B)> (1+\xi)n.$
It is clear that
 $$(1+\xi)n \le|V_1|\le (N-1)/2\le(2+\xi/2)n.$$
Since $\delta(B[V_1])\ge(1+\xi)n-1>|V_1|/2$, it follows from Lemma \ref{1-thm-1} that $B[V_1]$ is pancyclic,
which implies that there is a blue $C_n$.

Now, for any $0<\xi<1/10$, if we take $n_0=\lceil1000/(3\xi)\rceil$, then
$\delta(B)\geq (1+\xi)n\geq \frac{N}{4}+250$ for all $n\ge n_0$.
 Therefore, we can find a blue $C_n$ by Lemma \ref{1-thm-3}, Lemma \ref{1-thm-2} and Fact \ref{f1}.
The proof for $k=3$ is complete.\hfill$\Box$

\bigskip\noindent
{\large\bf Step 2: $B_n^{(k)}$ versus $C_n$ for $k\ge4$}
\medskip

For $k=3$, the result has been verified for  sufficiently large  $n$ in Step 1. We now suppose that the assertion holds for some $k\ge3$ and  prove  it  for $k+1$.

Let $0<\xi<1/10$ be fixed and $N=\lfloor(k+2)(1+\xi)n\rfloor$. We consider a red-blue edge coloring of $K_N$ on vertex set $V$.
In the following, we will omit the ceiling and floor  as it will not affect the result.
 We choose $\delta$ sufficiently small. To be precise,
 \begin{align}\label{dlt}
0<\delta<\min\left\{10^{-15},\;\frac{\xi^2}{600(k+2)^2}\right\}.
 \end{align}
Lemma \ref{con} with $\delta$ and $k+1$ gives us a constant $\eta$.
Let $\beta$ and $\epsilon$ be  sufficiently small such that
 \begin{align}\label{beta-epi}
\beta=\frac{\xi}{20(k+2)^2} \;\;\mbox{and}\;\;0<\epsilon
<\min\left\{\eta,\;\beta^2\right\}.
 \end{align}
  Set
 \begin{align}\label{alpha}
 \alpha=\frac{1}{k+2}-\beta-\sqrt{\epsilon}.
 \end{align}
Let $M$ be given by  Lemma \ref{con0} with $\epsilon$ and large $m_0=\lceil\frac{1}{\epsilon}\rceil$.
We apply Lemma \ref{con0} to the red subgraph $R$ and obtain an equitable partition $V (G) = \sqcup^m_{i=1}V_i$ and subsets $W_i \subset V_i$ such that $W_i$ is $\epsilon$-regular for all $i$ and, for all but $\epsilon m^2$ pairs $(i,j)$ with $1\leq i \neq j \leq m$, $(V_i,V_j)$ and $(W_i, V_j)$ are $\epsilon$-regular with $|d_R(W_i,V_j)-d_R(V_i,V_j)| \leq \epsilon$.
Here we do not require pairs $(W_i,W_j)$ for $1\le i,j\le m$ to be $\epsilon$-regular.
For  convenience, we will assume $|V_i|=\frac{N}{m}$ for all $1 \leq i \leq m$. If $n$ is large enough, then $\frac{N}{m} \geq \frac{N}{M} \geq \max\{n_1,n_2\}$, where $n_1$ is given by Lemma \ref{ben-s} with $\beta-\epsilon$, and $n_2$ is given by Lemma \ref{tluczak} with $\delta$ and $\alpha$.
 Note that a partition obtained by applying Theorem \ref{con0} to $R$ is also such a partition for $B$.

Let $F$ be the reduced graph defined on $\{v_1,v_2,\dots,v_m\}$, in which $v_i$ and $v_j$ are non-adjacent in $F$ if the pairs $(V_i,V_j)$ and $(W_i, V_j)$ are not all $\epsilon$-regular with $|d_R(W_i,V_j)-d_R(V_i,V_j)| \leq \epsilon$. Then the number of edges of $F$ is at least
$(1-\epsilon)m^2$. Therefore, by deleting at most $\sqrt{\epsilon}m$ vertices, we may assume that each vertex is adjacent to at least $(1-\sqrt{\epsilon})m$ vertices.  In what follows, when referring
to the reduced graph, we will assume that these vertices have been removed. For each remaining vertex $v_i$, we color $v_i$ red if the density of the red subgraph induced by $W_i$ satisfies that $d_R(W_i)\ge1/2$, and we color $v_i$ blue otherwise.
We color an edge $v_iv_j$ red if $d_R(V_i,V_j)\ge 1-\beta$, or blue if $d_B(V_i,V_j)\ge \beta$.
  Let $F_R$ and $F_B$ be the subgraphs formed by red edges and blue edges of $F$, respectively.

For a blue  vertex $v_a$ in $F$, suppose $d_{F_R}(v_a)\ge m_1:=(\frac{k+1}{k+2}+\beta)m$.
Recall an edge $v_iv_j$ in $F$ is red if and only if $d_R(V_i,V_j)\ge 1-\beta$. By averaging, there exists a vertex $u\in V_a$ such that
\[
d_R(u)\ge (1-\beta)\frac{N}{m}\cdot m_1=(1-\beta)(k+2)(1+\xi)\left(\frac{k+1}{k+2}+\beta\right)n>(k+1)(1+\xi)n,
\]
here we used the assumption (\ref{beta-epi}).
Thus, by the induction hypothesis, there is either a blue $C_n$ or a red $B_n^{(k)}$ in $N_R(u)$. In the former case, we are through. In the latter case,
the red $B_n^{(k)}$ together with $u$ gives a red $B_n^{(k+1)}$. We are done with this case.
Therefore, we assume that each blue vertex $v_a$ in the reduced graph $F$ satisfies
\begin{align}\label{blue-vt}
d_{F_B}(v_a)\ge m_2:=\left(\frac{1}{k+2}-\beta-\sqrt{\epsilon}\right)m.
\end{align}

For a red vertex  $v_a$ in $F$, if $d_{F_R}(v_a) \geq \frac{(1-0.5\xi)m}{k+2}$, then we can apply Lemma \ref{con} to find a red $B_n^{(k+1)}$.  Recall  the red density $d_R(W_a) \geq 1/2>\delta$, where $W_a\subseteq V_a$. We apply Lemma \ref{con} with $U_i=W_a$ for $1 \leq i \leq k+1$ and $U_{k+1+s}$ equal to each of the $V_s$ for which $(W_a,V_s)$ is $\epsilon$-regular and the edge $v_av_s$ is red. We conclude that there is a red $K_{k+1}$ which is contained in at least
\begin{align*}
\sum_{s} \left([d_R(W_a,V_s)]^{k+1}-\delta\right) |V_s|
&\geq \left((1-\beta-\epsilon)^{k+1}-\delta\right)\frac{N}{m} \cdot \frac{(1-0.5\xi)}{k+2}m
\\&\ge\left((1-\beta-\epsilon)^{k+1}-\delta\right)\left(1+\frac{\xi}{3}\right) n
\\&\ge\left(1+\frac{\xi}{3}-2(k+1)(\beta+\epsilon)-2\delta\right) n
\end{align*}
red $K_{k+2}$. We notice this quantity is at least $n$ since $\beta,\epsilon$ and $\delta$ are sufficiently small in terms of $1/k$ and $\xi$  from (\ref{dlt}) and (\ref{beta-epi}), so we are through as there is a red $B_n^{(k+1)}.$ Therefore, for the rest of the proof, we assume that each red vertex $v_a$ satisfies
\begin{align}\label{red-vt}
d_{F_B}(v_a) \geq \left(\frac{k+1+0.5\xi}{k+2}-\sqrt{\epsilon}\right)m.
\end{align}

From (\ref{blue-vt}) and (\ref{red-vt}), we have $d_{F_B}(v_a)\ge m_2:=(\frac{1}{k+2}-\beta-\sqrt{\epsilon})m$
 no matter which color a vertex $v_a$ has received.
We separate the proof into two cases depending on the parity of $n$.

\medskip
\noindent
{\bf Case 1.} $n$ is even.
 \medskip

We apply the techniques used by {\L}uczak \cite{Luczak-1999} and Lemma \ref{ben-s} to show that the blue graph $B$ contains a $C_n$.
By Erd\H{o}s--Gallai theorem \cite{EG}, $F_B$ contains  a path  $v_1v_2\cdots v_{\ell}$, here $\ell\ge m_2$.
For each $1 \leq i \leq \ell$, we split  $V_i$ corresponding to the vertex $v_i$ into two subsets $V_i'$ and $V_i''$,
where $|V_i'|,|V_i''|\ge |V_i|/2= \frac{N}{2m}$.
We have an even ``fat'' cycle $V_1'V_2'\cdots V_{\ell}'V_{\ell-1}''V_{\ell-2}''\cdots V_2''V_1'$.
For convenience, we relabel it as $U_1U_2\cdots U_{\ell}U_{\ell+1}U_{\ell+2}\cdots U_{2\ell-2}U_1$.
Note that Fact \ref{hdt} implies that $(U_i,U_{i+1})$ is $3\epsilon$-regular with blue density at least $\beta-\epsilon$ for $1\le i\le2\ell-2$, where the sums of the indices are taken modulo $2\ell-2$.
By Fact \ref{lar-dg}, there are at least $(1-6\epsilon)|U_i|$ vertices of $U_i$ having $(\beta-4\epsilon) \tfrac{N}{2m}$ neighbors in each sets of $U_{i-1}$ and $U_{i+1}$. Therefore, for each $1 \leq i \leq 2\ell-2$, we can choose a vertex $u_i \in U_i$ such that $u_1u_2\cdots u_{2\ell-2}u_1$ form an even cycle satisfying $d_B(u_i,U_{i-1}) \geq 4\beta/5 \cdot \tfrac{N}{2m}$ and $d_B(u_i,U_{i+1}) \geq 4\beta/5 \cdot \tfrac{N}{2m}$.
By Lemma \ref{ben-s}, for each $1 \leq i \leq \ell-1$, we can replace the edge $u_{2i-1}u_{2i}$ by an odd path with  endpoints $u_{2i-1}$ and $u_{2i}$ using vertices from $U_{2i-1}$ and $U_{2i}$, here the length of the path can vary from 1 to $(1-5\epsilon/\beta)\cdot \tfrac{N}{2m}$. Therefore, we can enlarge this  even cycle to all even cycles of length from $2\ell-2$ to $(2\ell-2)(1-5\epsilon/\beta)\frac{N}{2m}$.
 Recall  assumptions (\ref{beta-epi}), $m\ge m_0=\lceil\frac{1}{\epsilon}\rceil$ and $\ell\ge m_2=\left(\frac{1}{k+2}-\beta-\sqrt{\epsilon}\right)m$.  We have  $(2\ell-2)(1-5\epsilon/\beta)\frac{N}{2m}\ge n$. Therefore, there is a blue $C_n$ as desired.

\medskip
\noindent
{\bf Case 2.} $n$ is odd.
 \medskip

First, we suppose that $F_B$ contains an odd cycle  $v_1v_2\ldots v_{2\ell+1}v_1$ such that $2\ell+1\ge m_2$.
 Similarly, by Fact \ref{hdt}, we can find an odd cycle $u_1u_2\ldots u_{2\ell+1}u_1$ such that $u_i \in V_i$ , $d_{B}(u_i,V_{i-1}) \geq  4\beta/5 \cdot |V_{i-1}|$, and $d_{B}(u_i,V_{i+1}) \geq 4\beta/5 \cdot |V_{i+1}|$ for each $1 \leq i \leq 2\ell+1$.
 Now, applying Lemma \ref{ben-s}, for each $1 \leq i \leq \ell$, we replace the edge  $u_{2i-1}u_{2i}$ by odd path using vertices from $V_{i-1}$ and $V_i$ to enlarge this odd cycle to  odd cycles of length from $2\ell+1$ to $2\ell(1-5\epsilon/\beta)\frac{N}{m}$ in the blue graph $B$.
We can definitely find an odd cycle $C_n$ as $2\ell(1-5\epsilon/\beta)\frac{N}{m} \geq n$.

It remains to consider the case where  $F_B$ contains no odd cycles of length at least $m_2$. Erd\H{o}s--Gallai theorem \cite{EG} already implies that $F_B$ contains an even cycle of length at least $m_2$. We assume that $C:=v_1v_2,\cdots v_{2s} v_1$ is such an even cycle, here $2s \ge m_2$.

\medskip
{\bf Subcase 2.1.} There exists a blue vertex in $C$, say $v_{2s}$.
 \medskip

 Recall the definition of $W_{2s}$. We note $(V_{2s-1},W_{2s})$ and $(W_{2s},V_1)$ are $\epsilon$-regular pairs. Moreover, $d_B(V_{2s-1},W_{2s}) \geq \beta-\epsilon$ and $d_B(W_{2s},V_{1}) \geq \beta-\epsilon$.  We apply Lemma \ref{edge} with $A=V_{2s-1}, B=W_{2s},$ $C=V_{1}$ and $\beta'=\beta-\epsilon$ to get a blue edge $xy$ in $W_{2s}$ which satisfies $d_B(x,V_{2s-1}) \geq (\beta'-\epsilon) \cdot\tfrac{N}{m} \geq 4\beta/5 \cdot \tfrac{N}{m}$ and
  $d_B(y,V_{1}) \geq (\beta'-\epsilon) \tfrac{N}{m}$.  Similar to the case where $F_B$ contains an odd blue cycle of length at least $m_2$, we can find a blue odd cycle $u_1 u_2 u_3 \ldots u_{2s+1} u_1$, where $u_i \in V_i$ for $1 \leq i \leq 2s-1$, $u_{2s}=x$, and $u_{2s+1}=y$. Furthermore, we have
$d_B(u_1,V_2) \geq   4\beta/5 \cdot \tfrac{N}{m}$. For each $2 \leq i \leq 2s-1$, it satisfies $d_B(u_{i},V_{i-1}) \geq 4\beta/5 \cdot \tfrac{N}{m}$ and $d_B(u_{i},V_{i+1}) \geq 4\beta/5 \cdot \tfrac{N}{m}$. For each $1 \leq i \leq s$,  we again apply Lemma \ref{ben-s} to replace each blue edge $u_{2i-1}u_{2i}$ by an odd blue path. When we enlarge the edge $u_{2s-1}u_{2s}$, we use vertices from $V_{2s-1} \cup (V_{2s}\setminus \{y\})$ to avoid the vertex $u_{2s+1}=y \in V_{2s}$.
Thus $F_B$ contains odd cycles of length from $2s+1$ to  $2s(1-5\epsilon/\beta)\tfrac{N}{m}$ and there is a blue $C_n$ as $2s(1-5\epsilon/\beta)\tfrac{N}{m} \geq n$.

\medskip
{\bf Subcase 2.2.} Each vertex of $C$ is red.
 \medskip

Write $d=2s$ for the length of the cycle $C$.
   Recalling (\ref{red-vt}), we have a lower bound on the number of edges in $F_B$ as follows:
\begin{align*}
|E(F_B)| \geq\sum_{i=1}^d d_{F_B}(v_i) -\binom{d}{2}
\geq \left(\frac{k+1+0.5\xi}{k+2}-\sqrt{\epsilon}\right)dm-\binom{d}{2}.
\end{align*}
Since $d \geq m_2$ and the right hand side above is an increasing function of $d$ when $d \geq m_2$,
 it follows that
\begin{align} \label{lbedge}
|E(F_B)| &\geq \nonumber\left(\frac{k+1+0.5\xi}{k+2}-\sqrt{\epsilon}\right)\left(\frac{1}{k+2}-\beta-\sqrt{\epsilon}\right)m^2
-\frac{1}{2}\left(\frac{1}{k+2}-\beta-\sqrt{\epsilon}\right)^2m^2
\\& \nonumber \geq \left(\frac{2k+1}{2(k+2)^2}+\frac{0.5\xi}{(k+2)^2} -\beta-\sqrt{\epsilon}- \frac{\sqrt{\epsilon}}{k+2}\right)m^2
\\&\geq \left(\frac{2k+1}{2(k+2)^2}+\frac{\xi}{3(k+2)^2}  \right)m^2,
\end{align}
 here we note that $\beta$ and $\sqrt{\epsilon}$ are sufficiently small in terms of $1/k^2$ and $\xi$ from (\ref{beta-epi}).

We next apply Lemma \ref{tluczak} to obtain an upper bound on $|E(F_B)|$.  Actually, we will apply Lemma \ref{tluczak} with $G=F_B$, $\delta$ and $\alpha=\frac{1}{k+2}-\beta-\sqrt{\epsilon}$. Since $F_B$ contains no odd cycles of length at least $m_2=\alpha m$, there are two subgraphs $G'=F'_B$ and $G''=F''_B$ satisfying all properties listed in the lemma. If $F'_B$ is empty, then we get $|E(F_B)| \leq \frac{1}{2}(\frac{1}{k+2}-\beta-\sqrt{\epsilon})m^2+\delta m^2$, which clearly is a contradiction to the lower bound on $|E(F_B)|$ from (\ref{lbedge}). If $F_B'$ is not empty, then we assume that $X$ and $Y$ are two color classes of $F_B'$.


\medskip
{\em Claim.} We have $\max\{|X|,|Y|\} \leq k+1+\frac{(1-0.5\xi)m}{k+2}.$
 \medskip

  {\em Proof.}
  Suppose  $|X| \geq k+1+\frac{(1-0.5\xi)m}{k+2}$ without loss of generality.
  We aim to find a book $B_{n}^{(k+1)}$ in the red graph.
Since all but at most $\delta m^2$ edges of $F_B$ belong to either $F'_B$ or $F''_B$  by (4) of Lemma \ref{tluczak}, it follows that  $X$ contains at least
${|X|\choose2}-(\epsilon +\delta) m^2\ge (1-\gamma){|X|\choose2}$ red edges, where
 $0\le\gamma<3(k+2)^2(\epsilon +\delta)$.
  By deleting at most $\sqrt{\gamma}|X|$ vertices, we may assume that each vertex of $X$ is red-adjacent to at least $(1-\sqrt{\gamma})|X|$ vertices in $X$. We notice
  $$(1-(k+1)\sqrt{\gamma})|X|\ge m_3:=\frac{(1-0.6\xi)m}{k+2},$$
  here we note $\epsilon,\delta<\frac{\xi^2}{600(k+2)^2}$ from (\ref{dlt}) and (\ref{beta-epi}).
 Therefore, we can choose vertices $v_1,\dots,v_{k+1}$ from $X$ step by step such that
  $v_1,\dots,v_{k+1}$ form a red clique in the reduced graph $F$ and
$v_1,\dots,v_{k+1}$ have at least $m_3$ common red neighbors, say $v_{k+2},\dots,v_{k+1+m_3}$. i.e., $X$ contains a red book $B_{m_3}^{(k+1)}$.
We apply Lemma \ref{con} with $U_i=V_i$ for $1 \leq i \leq k+1$ and $U_{k+1+j}=V_{k+1+j}$ for $1\le j\le m_3$,
and conclude that there is a red $K_{k+1}$ which is contained in at least
\begin{align*}
\sum_{j=1}^{m_3} \left(\prod_{i=1}^kd_R(V_i,V_{k+1+j})-\delta\right) |V_{k+1+j}|
&\geq \left((1-\beta)^{k+1}-\delta\right)\frac{N}{m} \cdot \frac{(1-0.6\xi)}{k+2}m
\\&\ge\left((1-\beta)^{k+1}-\delta\right)\left(1+\frac{\xi}{3}\right) n
\end{align*}
red $K_{k+2}$.
This quantity is at least $n$ since $\beta$ and $\delta$ are sufficiently small in terms of $1/k$ and $\xi$ from (\ref{dlt}) and (\ref{beta-epi}).  We found a red $B_n^{(k+1)}$ as desired. \hfill$\Box$

\medskip
Consequently, $|E(F_B')| \leq \big( k+1+\frac{(1-0.5\xi)m}{k+2}\big)^2$ from the above claim. The number of edges in $F_B$ can be bounded from above as follows:
\begin{align*}\label{ubedge}
|E(F_B)| &\leq \nonumber \left(k+ 1+\frac{(1-0.5\xi)m}{k+2}\right)^2+\frac{1}{2}\left(\frac{1}{k+2}-\beta-\sqrt{\epsilon}\right)m^2+\delta m^2
\\&\nonumber\le\left(\frac{k+4}{2(k+2)^2}-\frac{\xi}{(k+2)^2}+\delta\right)m^2+2m+(k+1)^2
\\&\le\left(\frac{k+4}{2(k+2)^2}-\frac{\xi}{2(k+2)^2}\right)m^2,
\end{align*}
 here we note that $m\ge m_0\ge 1/\epsilon$,   $\epsilon$ and $\delta$ are sufficiently small in terms of $1/k^2$ and $\xi$ from (\ref{dlt}) and (\ref{beta-epi}).
Recalling the inequality \eqref{lbedge}, we obtained
 \begin{equation}\label{contra}
 \left(\frac{2k+1}{2(k+2)^2}+\frac{\xi}{3(k+2)^2}  \right)m^2 \leq |E(F_B)| \leq \left(\frac{k+4}{2(k+2)^2}-\frac{\xi}{2(k+2)^2}\right)m^2.
 \end{equation}
This is a contradiction provided $k \geq 3$.

This completes the proof of the induction step and hence Theorem \ref{main}.
\hfill$\Box$

\noindent
{\bf Remark 2:} For $k=2$, the inequality \eqref{contra} indeed holds  and there is no contradiction.
Thus we are not able to prove $r(B_n^{(3)},C_n)=(4+o(1))n$ by induction and we provide a separated proof for this case.

\noindent
{\bf Acknowledgement:} The authors would like to thank anonymous referees for their valuable comments which greatly improve the presentation of this paper.

\end{spacing}

\end{document}